\documentclass[11pt]{article}
\usepackage{amsfonts,mathrsfs,amssymb,amsthm,mathptm}
\usepackage{amsmath,amscd}
\usepackage{mathptm,pslatex}
\usepackage{multirow}
\usepackage{color}
\usepackage[dvips]{graphicx}
\usepackage[all]{xy}
\usepackage{graphicx}
\usepackage{fancyhdr}
\usepackage{url}

\begin{document}
\newtheorem{Def}{Definition}[section]
\newtheorem{Bsp}[Def]{Example}
\newtheorem{Prop}[Def]{Proposition}
\newtheorem{Theo}[Def]{Theorem}
\newtheorem{Lem}[Def]{Lemma}
\newtheorem{Koro}[Def]{Corollary}
\theoremstyle{definition}
\newtheorem{Rem}[Def]{Remark}

\newcommand{\add}{{\rm add}}
\newcommand{\con}{{\rm con}}
\newcommand{\gd}{{\rm gl.dim}}
\newcommand{\dm}{{\rm domdim}}
\newcommand{\tdim}{{\rm dim}}
\newcommand{\E}{{\rm E}}
\newcommand{\Mor}{{\rm Morph}}
\newcommand{\End}{{\rm End}}
\newcommand{\ind}{{\rm ind}}
\newcommand{\rsd}{{\rm res.dim}}
\newcommand{\rd} {{\rm rep.dim}}
\newcommand{\ol}{\overline}
\newcommand{\overpr}{$\hfill\square$}
\newcommand{\rad}{{\rm rad}}
\newcommand{\soc}{{\rm soc}}
\renewcommand{\top}{{\rm top}}
\newcommand{\stp}{{\mbox{\rm -stp}}}
\newcommand{\pd}{{\rm projdim}}
\newcommand{\id}{{\rm injdim}}
\newcommand{\fld}{{\rm flatdim}}
\newcommand{\fdd}{{\rm fdomdim}}
\newcommand{\Fac}{{\rm Fac}}
\newcommand{\Gen}{{\rm Gen}}
\newcommand{\fd} {{\rm findim}}
\newcommand{\Fd} {{\rm Findim}}
\newcommand{\Pf}[1]{{\mathscr P}^{<\infty}(#1)}
\newcommand{\DTr}{{\rm DTr}}
\newcommand{\cpx}[1]{#1^{\bullet}}
\newcommand{\D}[1]{{\mathscr D}(#1)}
\newcommand{\Dz}[1]{{\mathscr D}^+(#1)}
\newcommand{\Df}[1]{{\mathscr D}^-(#1)}
\newcommand{\Db}[1]{{\mathscr D}^b(#1)}
\newcommand{\C}[1]{{\mathscr C}(#1)}
\newcommand{\Cz}[1]{{\mathscr C}^+(#1)}
\newcommand{\Cf}[1]{{\mathscr C}^-(#1)}
\newcommand{\Cb}[1]{{\mathscr C}^b(#1)}
\newcommand{\Dc}[1]{{\mathscr D}^c(#1)}
\newcommand{\K}[1]{{\mathscr K}(#1)}
\newcommand{\Kz}[1]{{\mathscr K}^+(#1)}
\newcommand{\Kf}[1]{{\mathscr  K}^-(#1)}
\newcommand{\Kb}[1]{{\mathscr K}^b(#1)}
\newcommand{\modcat}{\ensuremath{\mbox{{\rm -mod}}}}
\newcommand{\Modcat}{\ensuremath{\mbox{{\rm -Mod}}}}

\newcommand{\stmodcat}[1]{#1\mbox{{\rm -{\underline{mod}}}}}
\newcommand{\pmodcat}[1]{#1\mbox{{\rm -proj}}}
\newcommand{\imodcat}[1]{#1\mbox{{\rm -inj}}}
\newcommand{\Pmodcat}[1]{#1\mbox{{\rm -Proj}}}
\newcommand{\Imodcat}[1]{#1\mbox{{\rm -Inj}}}
\newcommand{\opp}{^{\rm op}}
\newcommand{\otimesL}{\otimes^{\rm\mathbb L}}
\newcommand{\rHom}{{\rm\mathbb R}{\rm Hom}\,}
\newcommand{\projdim}{\pd}
\newcommand{\Hom}{{\rm Hom}}
\newcommand{\Coker}{{\rm Coker}}
\newcommand{ \Ker  }{{\rm Ker}}
\newcommand{ \Cone }{{\rm Con}}
\newcommand{ \Img  }{{\rm Im}}
\newcommand{\Ext}{{\rm Ext}}
\newcommand{\StHom}{{\rm \underline{Hom}}}

\newcommand{\gm}{{\rm _{\Gamma_M}}}
\newcommand{\gmr}{{\rm _{\Gamma_M^R}}}

\def\vez{\varepsilon}\def\bz{\bigoplus}  \def\sz {\oplus}
\def\epa{\xrightarrow} \def\inja{\hookrightarrow}

\newcommand{\lra}{\longrightarrow}
\newcommand{\llra}{\longleftarrow}
\newcommand{\lraf}[1]{\stackrel{#1}{\lra}}
\newcommand{\llaf}[1]{\stackrel{#1}{\llra}}
\newcommand{\ra}{\rightarrow}
\newcommand{\dk}{{\rm dim_{_{k}}}}

\newcommand{\colim}{{\rm colim\, }}
\newcommand{\limt}{{\rm lim\, }}
\newcommand{\Add}{{\rm Add }}
\newcommand{\Tor}{{\rm Tor}}
\newcommand{\Cogen}{{\rm Cogen}}
\newcommand{\Tria}{{\rm Tria}}
\newcommand{\tria}{{\rm tria}}

{\Large \bf
\begin{center}
Cellularity of centrosymmetric matrix algebras and Frobenius extensions
\end{center}}

\medskip
\centerline{\textbf{Changchang Xi$^*$ and Shujun Yin}}

\renewcommand{\thefootnote}{\alph{footnote}}
\setcounter{footnote}{-1} \footnote{$^*$Corresponding author's
Email: xicc@cnu.edu.cn; Fax: 0086 10 68903637.}
\renewcommand{\thefootnote}{\alph{footnote}}
\setcounter{footnote}{-1} 
\renewcommand{\thefootnote}{\alph{footnote}}
\setcounter{footnote}{-1} 

\begin{abstract}
Centrosymmetric matrices of order $n$ over an arbitrary algebra $R$ form a subalgebra of the full $n\times n$ matrix algebra over $R$. It is called the centrosymmetric matrix algebra of order $n$ over $R$ and denoted by $S_n(R)$. We prove (1) $S_n(R)$ is Morita equivalent to $S_2(R)$ if $n$ is even, and to $S_3(R)$ if $n\ge 3$ is odd; (2) the full $n\times n$ matrix algebra over $R$ is a separable Frobenius extension of $S_n(R)$; and (3) if $R$ is a commutative ring, then $S_n(R)$ is a cellular $R$-algebra in the sense of Graham-Lehrer for all $n\ge 1$.
\end{abstract}

\section{Introduction}
Given a natural number $n$ and an associative (not necessarily commutative) ring $R$ with identity,
all symmetric sequences of $n^2$
elements in $R$ form a subring of the full $n\times n$ matrix ring $M_n(R)$. This subring is called the centrosymmetric matrix ring of $R$, denoted by $S_n(R)$.  The study of centrosymmentric matrices over a field has a long history in matrix theory (see \cite[p.19]{Muir}, \cite[p.124]{A}). A short survey of this topic in earlier time can be found in \cite{W} where, in particular, the structure of $\mathbb{R}$-algebra $S_n(\mathbb{R})$ was determined, that is, $S_n(\mathbb{R})\simeq M_{n-k}(\mathbb{R})\times M_{k}(\mathbb{R})$, with $k=\lceil n/2\rceil$, the least integer bigger than or equal to $n/2$. Thus $S_n(\mathbb{R})$ is always a semisimple $\mathbb{R}$-algebra. Though there are lots of works on different aspects of matrix theory of centrosymmetric matrices, less attentions are focus on algebraic aspects of such matrix algebras.

In this note, we shall investigate ring-theoretical structures of the centrosymmetric matrix algebras $S_n(R)$ for an arbitrary ring $R$. This will be carried out by using methods in both matrix theory and representation theory of algebras. It turns out that $S_n(R)$ may not be semisimple even if $R$ is a field. We also study relations between $S_n(R)$ and $M_n(R)$.  Our main result is the following.

\medskip
{\bf Theorem.} (1) {\it For any ring $R$, $M_n(R)$ is a separable Frobenius extension of $S_n(R)$; and $S_n(R)$ is Morita equivalent to $S_2(R)$ if $n$ is even, and to $S_3(R)$ if $n\ge 3$ is odd.}

(2) {\it If $R$ is a commutative ring, then $S_n(R)$ is a cellular $R$-algebra for $n\ge 1$.}

\medskip
We mention that Frobenius extensions have connections with many aspects in mathematics from classical Yang-Baxter equation, Jones polynomials, Hecke algebras, quantised enveloping algebras to topological quantum field theories (see \cite{BGS, Kad, Kock2004}).

The proof of this result is given in Section \ref{sect3} after some preparations in Section \ref{sect2} where we give the precise definitions of all terminology unexplained in the result and all elementary facts needed in the proof.

\section{Preliminaries\label{sect2}}

Throughout this note, $R$ stands for a ring with identity. We denote by $M_n(R)$ the ring of all $n\times n$ matrices over $R$ for any positive integer $n$, and by $R\Modcat$ the category of all left $R$-modules. Recall that two rings $R$ and $S$ are \emph{Morita equivalent} if the categories $R\Modcat$ and $S\Modcat$ are equivalent, denoted by $R\approx S$. For example, $R$ is Morita equivalent to $M_n(R)$ for all $n\ge 1$.

A sequence $a=(a_1,\cdots, a_{n})$ of $n$ elements $a_1, a_2, \cdots, a_{n}$ in $R$ is said to be \emph{symmetric} if $a_i=a_{n+1-i}$ for all $1\le i\le n$. The number $n$ is called the length of the sequence $a$. Let $\mathcal{S}_n(R)$ denote the set of all symmetric sequences of length $n$ in $R$.

Now we define
$$S_n(R):=\{a\in M_n(R)\mid a_{i j}=a_{n+1-i,n+1-j}, 1\le i,j\le n\}.$$
Matrices in $S_n(R)$ are called \emph{centrosymmetric} matrices. Examples of centrosymmetric matrices include symmetric Toeplitz matrices.

If $m=n^2$, then each element $a\in \mathcal{S}_m(R)$ gives rise to a unique $n\times n$ matrix $A:=(a_{ij})\in S_n(R)$ with the first row $(a_1, \cdots, a_n)$, the second row $(a_{n+1},\cdots, a_{2n}), \cdots,$ and the $n$-th row $(a_{(n-1)n+1}, \cdots, a_{n^2})$. So we identify $a$ with $A$ and consider $\mathcal{S}_m(R)$ as a subset of $M_n(R)$, and therefore $\mathcal{S}_m(R)$ can be identified with $S_n(R).$

Let $e_{ij}$ be the matrix units of $M_n(R)$ with $1\le i,j\le n$, and let $c:=e_{1 n}+e_{2,n-1}+\dots + e_{n-1,2}+e_{n 1}$ be the skew diagonal identity in $M_n(R)$. Then $c\in S_n(R)$ and $c^2=1$. Note that, for any $a=(a_{ij})\in M_n(R)$, we have
$$(cac)_{jk}=a_{n+1-j,n+1-k}, \; 1\le j,k\le n.$$
Recall that a matrix $a=(a_{ij})\in M_n(R)$ is called \emph{persymmetric} if $ca'c=a$ where $a'$ stands for the transpose of the matrix $a$, that is, $a_{ij}=a_{n+1-j, n+1-i}$ for all $1\le i,j\le n$.  Thus, if $a$ is \emph{bisymmetric}, that is, both symmetric and persymmetric, then $a\in S_n(R)$. Let $B_n(R)$ be the set of all bisymmetric matrices in $M_n(R)$. Then $B_n(R)\subseteq S_n(R)$. In general, $B_n(R)$ may not be closed under multiplication (for example, in $M_3(\mathbb{Z}), (e_{11}+e_{13}+e_{31}+e_{33})(e_{12}+e_{21}+e_{23}+e_{32})=2(e_{12}+e_{32})$ is not bisymmetric). But we will see that $S_n(R)$ is closed under multiplication.

Centrosymmetric and bisymmetric matrices over fields have been intensively studied in matrix theory. For information, we may refer the reader to \cite{A, BM, CB, Muir, W}.

\medskip
Now we define an anti-Kronecker delta $\sigma$: $\sigma_{ij}=0$ if $i=j$, and $\sigma_{ij}=1$ if $i\neq j$. Further, we define
$$ f_{ii}:= e_{ii}+\sigma_{i, n+1-i}e_{n+1-i, n+1-i}, \quad f_{ij}:=e_{ij}+e_{n+1-i,n+1-j}, 1\le i\neq j\le n.$$
Then $\{f_{ij}\mid 1\le i,j\le \lceil n/2 \rceil\}$ is an $R$-basis of $S_n(R)$ and $1=\sum_{i=1}^{\lceil n/2 \rceil} f_{ii}$ is a decomposition of $1$ into pairwise orthogonal idempotents in $S_n(R)$. In the sequel, we often write $f_i$ for $f_{ii}$ for simplicity.

\begin{Lem}\label{lem1} $(1)$ An element $a=(a_{ij})\in M_n(R)$ belongs to $S_n(R)$ if and only if $cac=a$.

$(2)$ $S_n(R)$ is a subring of the matrix ring $M_n(R)$.

$(3)$ $S_n(R)$ is a free $R$-module of rank $\lceil n^2/2 \rceil$, where $\lceil m\rceil$ denotes the least natural number $z$ such that $m\le z$.
\end{Lem}

{\it Proof.} (1) is obvious. (3) follows from the identification of $\mathcal{S}_{n^2}(R)$ with $S_n(R)$.

(2) This was observed for $ R=\mathbb{R}$ in \cite{W}. It is easy to check: if $a, b\in S_n(R)$, then $ab=(cac)(cbc)=cac^2bc=c(ab)c$, that is, $ab\in S_n(R).$ $\square$

\begin{Lem} \label{lem2}
 $(1)$ $ f_iS_n(R)f_j = Rf_{ij}+Rf_{i, n+1-j}, \quad f_{i, n+1-j}^2=\delta_{i j}f_{ij}+ \delta_{i, n+1-j}f_{i,n+1-j}.$

$(2)$ $ f_{ij}f_{pq}=\delta_{j p}f_{i q}+\delta_{j, n+1-p}f_{i, n+1-q},$ where $\delta_{ij}$ is the Kronecker delta.

$(3)$  $f'_{ij}=f_{ji}, \; \mbox{ and } f'_{ij}=f_{ij} \mbox{ if and only if either } i=j \mbox{ or  }  i+j=n+1.$
\end{Lem}

Thus the action of the transpose on $f_iS_n(R)f_i$  is identity.

\begin{Bsp}\label{bsp} $(1)$ {\rm $S_1(R)=R$ and $S_2(R)\simeq R[C_2]$, the group algebra of the cyclic group $C_2$ of order $2$ over $R$. In fact, $f_1f_{12}=f_{12}=f_{12}f_1$, $f_{12}^2=f_1$. Thus $S_2(R)=Rf_1+Rf_{12}\simeq R[C_2]$. Note that if $R$ is a field, then $R[C_2]\simeq R\times R$ if $R$ is not of characteristic $2$, and $R[C_2]\simeq R[X]/(X^2)$ if $R$ is of characteristic $2$.}\end{Bsp}

(2) $S_3(R)\simeq  \begin{pmatrix}R[C_2] & R\\ R & R\end{pmatrix}$, where the multiplication is given by
$$ \begin{pmatrix} a+bx & u \\ d & v\end{pmatrix} \begin{pmatrix} a_1+b_1x & u_1\\ d_1 & v_1\end{pmatrix}\begin{pmatrix} aa_1+bb_1+ ud_1+ (ab_1+ba_1+ud_1)x & au_1+bu_1+uv_1\\ da_1+db_1+vd_1 & 2du_1+vv_1\end{pmatrix}$$
and where $x$ is a generator of $C_2$.

In fact, $\{f_1, f_2, f_{12}, f_{13}, f_{21}\}$ is an $R$-basis of $S_3(R)$. An calculation shows that $f_1S_n(R)f_1=Rf_1+Rf_{13}
$ with $f_{13}^2=f_1$, $f_{12}f_{21}=f_1+f_{13}$, and $f_2S_3(R)f_2=Rf_2$, $f_{21}f_{12}=2f_2.$ Thus we get the above multiplication of $S_3(R)$.

\medskip
Now we recall the definitions of Frobenius extensions and cellular algebras.

Let $A$ and $B$ be rings with identity. If $B$ is a subring of $A$ with the same identity, then we say that $A$ (or $B\subseteq A$) is an extension of $B$. One of the interesting extensions is the class of Frobenius extensions introduced initially by Kasch (see \cite{kasch}). They play an important role in topological quantum field theories in dimension $2$ and even $3$ (see \cite{Kock2004}) and in representation theory and knot theory (see \cite{Kadison1996}, \cite{Kad}, \cite{Xi} and the references therein). Also, each Frobenius extension $B\subseteq A$ with $B$ a commutative ring provides us with a series of solution to classical Yang-Baxter equation (see \cite[Chapters 4 and 5]{Kad}).

\begin{Def} $(1)$ An extension $B\subseteq A$ of rings is called a Frobenius extension if $_BA$ is a finitely generated projective $B$-module and $\Hom_B(_BA,B)\simeq {}_AA_B$ as $A$-$B$-bimodules.

$(2)$ An extension $B\subseteq A$ is said to be split if the $B$-bimodule ${}_BB_B$ is a direct summand of ${}_BA_B$, and separable if the multiplication map $A\otimes_BA\mapsto A$ is a split surjective homomorphism of $A$-bimodule.
\end{Def}

Frobenius extensions have the following properties (see \cite[Theorem 1.2, p.3; Corollaries 2.16-17, p.15]{Kadison1996} for proofs).

\begin{Lem}\label{lem3} Let $B\subseteq A$ be an extension of rings.

$(1)$ The extension is a Frobenius extension if and only if there exists a $B$-$B$-bimodule homomorphism $E\in \Hom_{B-B}(_BA_B, {}_BB_B), x_i, y_i\in A$, $1\le i\le n,$ such that, for any $a\in A,$
$$\sum_{i=1}^n x_iE(y_ia) = a = \sum_i E(ax_i)y_i.$$ In this case, $(E, x_i,y_i)$ is called a Frobenius system of the extension.

$(2)$ Suppose that $B\subseteq A$ is a Frobenius extension with a Frobenius system $(E,x_i,y_i)$. Then

(i) $B\subseteq A$ is separable if and only if there exists $d\in C_A(B):=\{a\in A\mid ab=ba, \forall b\in B\}$ such that $\sum_{i=1}^nx_idy_i=1$,

(ii) $B\subseteq A$ is split if and only if there exists $d\in C_A(B)$ such that $E(d)=1$.
\end{Lem}

Finally, we recall the definition of quasi-hereditary and cellular algebras introduced by Cline-Parshall-Scott (see, for example, \cite{CPS}) and by Graham-Lehrer (see \cite{GL}), respectively. The former was arised in the study of highest weight categories in representation theory of semisimple complex Lie algebras and of reductive algebraic groups, and the latter was motivated by Kashdan-Lusztig's canonical basis of Hecke algebras. For cellular algebras, we will use an equivalent formulation in \cite{KX}.

Let $R$ be a commutative ring and $A$ be an $R$-algebra, that is, $A$ is a ring with identity and there is a ring homomorphism $R\ra A$ such that its image lies in the center of $A$. For example, any ring is a $\mathbb{Z}$-algebra.

\begin{Def} $(1)$ An ideal $AeA$ of $A$ generated by an idempotent element $e\in A$ is called a heredity ideal if $eAe\simeq R$, $Ae_{eAe}$ and $_{eAe}eA$ are finitely generated, free modules over $eAe$, and the multiplication map $Ae\otimes_ReA\ra AeA$ is injective.
The algebra $A$ is said to be quasi-hereditary if there is a complete set of pairwise orthogonal idempotents, $\{e_1, e_2, \cdots, e_n\}$ in $A$, such that, by defining $J_i:=A(e_1+\cdots +e_{i})A$, the subquotient $J_i/J_{i-1}$ of the chain
$$ 0=J_0\subseteq J_1\subset\cdots\subset J_i\subset  J_{n-1}\subset J_n=A$$
is a heredity ideal in $A/J_{i-1}$ for each $i$.

$(2)$ Let $i:A\ra A$ be an $R$-involution (that is, an anti-automorphism of the $R$-algebra $A$ of order $2$).  A two-sided ideal $J$ in A is called a cell ideal if and only if $i(J)=J$ and there exists a left ideal $\Delta\subseteq J$ such that $\Delta$ is finitely generated and free over $R$ and that there is an isomorphism of $A$-bimodules, say $\alpha: J\simeq \Delta\otimes_R i(\Delta)$, where $i(\Delta)$ is the image of $\Delta$ under $i$, such that the diagram commutes:
$$\xymatrix@M=0.6mm{
  J\ar[r]^{\alpha}\ar[d]_{i} & *+[r]{\Delta\otimes_Ri(\Delta)}\ar[d]^{x\otimes y\mapsto i(y)\otimes i(x)}\\
  *+[r]{J} \ar[r]^{\alpha} & *+[r]{\Delta\otimes_Ri(\Delta)}
}$$
The algebra $A$ with the involution $i$ is called a cellular algebra if and only if
there is an $R$-module decomposition $A=J'_1\oplus J'_2\oplus \cdots \oplus J'_n$ with $i(J'_j)=J'_j$ for all $j$ such that putting $J_j=\oplus_{p=1}^j J'_p$ gives a chain
of two-sided ideals of $A$ and for each $j$ the quotient $J'_j=J_j/J_{j-1}$ is a cell ideal in $A/J_{j-1}$ with respect to the involution induced from $i$ on the quotient.
\end{Def}

The following fact is well known in the literature (see, for example, \cite{CPS, DR} or \cite{KX, KX2012}).
\begin{Lem} \label{lem2.7} Suppose that $A$ is an $R$-algebra and $e^2=e$ is an idempotent in $A$ such that $eAe\simeq R$ and $Ae$ is a free $R$-module of finite rank.
If $i$ is an $R$-involution on $A$ such that $i(e)=e$ and the multiplication map $Ae\otimes_{eAe}eA\ra AeA$ is injective, then $AeA$ is a cell ideal in $A$.
\end{Lem}

Clearly, the full matrix $R$-algebra $M_n(R)$ is quasi-hereditary and cellular with respect to the matrix transpose. The group algebra $R[C_2]$ of $C_2$ over $R$ is cellular with respect to the identity map since $R(1-x)$, $x^2=1$, is a cell ideal in $R[C_2]$ and the chain $R(1-x)\subset R[C_2]$ is a cell chain. If $R$ is a field of characteristic $2$, then $R[C_2]$ is not quasi-hereditary.

We remark that the cellularity of $R$-algebras does not have to be preserved under Morita equivalences.

\section{Proof of the result\label{sect3}}

We first prove the following result.

\begin{Theo} \label{thm1} Let $R$ be a ring. Then

$(1)$ $S_{2m+1}(R)$ is Morita equivalent to $S_3(R)$ for all $m\ge 1$.

$(2)$ $S_{2m}(R)$ is Morita equivalent to $R[C_2]$ for all $m\ge 1$, where $C_2$ is the cyclic group of order $2$.

$(3)$ $S_n(R)\subseteq M_n(R)$ is a separable Frobenius extension for all $n\ge 1.$

\end{Theo}

{\it Proof.} We assume $n\ge 4$ and define
$$\varphi: S_n(R)f_1\lra S_n(R)f_j, \quad af_1\mapsto af_1(e_{1 j}+e_{n,n+1-j})$$
for $1\le j\le \lfloor n/2\rfloor$, where $\lfloor p\rfloor$ denotes the largest natural number $z$ such that $z\le p$. Note that $f_1(e_{1 j}+e_{n,n+1-j})f_j=e_{1j}+e_{n,n+1-j}$. Thus $\varphi$ is well defined. Clearly, it is a homomorphism of left $S_n(R)$-modules.

Note that elements in $S_n(R)f_1$ and $S_n(R)f_j$ are of the forms
$$\begin{pmatrix}
a_{1 1} & 0 & \dots & 0 & a_{n 1} \\
a_{2 1} & 0 & \dots & 0 & a_{n-1,1} \\
\vdots & \vdots &  & \vdots & \vdots \\
a_{j 1} & 0 & \dots & 0 & a_{n+1-j,1} \\
\vdots & \vdots &  & \vdots & \vdots \\
a_{n 1} & 0 & \dots & 0 & a_{1 1} \\
\end{pmatrix}_{n\times n} \mbox{and}
\quad
\begin{pmatrix}
0& \dots & a_{1 j} & 0 & \dots & 0 & a_{n j}& \dots & 0 \\
0& \dots & a_{2 j} & 0 & \dots & 0 & a_{n-1,j} &\dots & 0 \\
0& \dots &\vdots & \vdots &  & \vdots & \vdots  &\dots & 0\\
0& \dots &a_{j j} & 0 & \dots & 0 & a_{n+1-j,j}& \dots & 0 \\
0& \dots &\vdots & \vdots &  & \vdots & \vdots & \dots & 0\\
0& \dots &a_{n j} & 0 & \dots & 0 & a_{1 j}& \dots & 0\\
\end{pmatrix}_{n\times n} (a_{i,j}\in R)$$
respectively.

Since $\varphi$ sends the first and $n$-th columns of $af_1$ to the $j$-th and $(n+1-j)$-th colums of $af_1(e_{1 j}+e_{n,n+1-j})$, respectively, it is clear that $\varphi$ is injective. Given an element $x\in S_n(R)f_j$, we define an element $a\in S_n(R)f_1$ by putting the $j$-th and $(n+1-j)$-th columns of $x$  as the first and $n$-th columns of $a$, and $0$ for other columns of $a$. Then $a$ is sent to $x$ by $\varphi$. Thus $\varphi$ is surjective.

(1) If $n=2m+1$, then $1=f_1+\cdots +f_m+ f_{m+1}$ in $S_{2m+1}(R)$. So, the above forms of matrices show $S_n(R)f_j\simeq S_n(R)f_1$ as left $S_n(R)$-modules for $2\le j \le m$. Hence $S_n(R)$ is Morita equivalent to $\End_{S_n(R)}\big(S_n(R)f_1\oplus S_n(R)f_{m+1}\big)$.

We shall prove $\End_{S_n(R)}\big(S_n(R)f_1\oplus S_n(R)f_{m+1}\big)\simeq S_3(R)$. Indeed, by definition, $f_{m+1}=e_{m+1, m+1}$. Clearly,

$$f_1S_n(R)f_1=\left\{\begin{pmatrix}
a_{11} & 0 & \dots & 0 & a_{1 n} \\
0 & 0 & \dots & 0 & 0 \\
\vdots & \vdots &  & \vdots & \vdots \\
0 & 0 & \dots & 0 & 0 \\
a_{1 n} & 0 & \dots & 0 & a_{11} \\
\end{pmatrix}\mid a_{ij}\in R \right\}=Rf_1+Rf_{1 n}\simeq R[C_2],$$

$$f_1S_n(R)f_{m+1}=\left\{\begin{pmatrix}
0     & \dots &0      & a_{1,m+1} &0       & \dots  & 0  \\
0     & \dots &0      & 0    &0   & \dots &0 \\
\vdots&  &\vdots  &\vdots     & \vdots &  & \vdots \\
0     & \dots &0      & 0    &0   & \dots &0 \\
0     & \dots &0      &a_{1,m+1}  & 0      & \dots  & 0\\
\end{pmatrix}\mid a_{1,m+1}\in R\right\}\simeq Rf_{1,m+1}, $$

$$f_{m+1}S_n(R)f_1=\left\{\begin{pmatrix}
0 & 0 & \dots & 0 & 0 \\
\vdots & \vdots &  & \vdots & \vdots \\
0 & 0 & \dots & 0 & 0 \\
a_{m+1,1} & 0 & \dots & 0 & a_{m+1,1} \\
0 & 0 & \dots & 0 & 0 \\
\vdots & \vdots &  & \vdots & \vdots \\
0 & 0 & \dots & 0 & 0 \\
\end{pmatrix}\mid a_{m+1,1}\in R\right\}\simeq Rf_{m+1,1}. $$
Further, we have
$$ f_{1,m+1}\in f_1S_n(R)f_{m+1}, \; f_{m+1,1}\in f_{m+1}S_n(R)f_1,\; f_{1,m+1}f_{m+1,1}=f_1+f_{1,n}, \; f_{m+1,1}f_{1, m+1}=2f_{m+1}.$$
Thus, by Example \ref{bsp}, $\End_{S_n(R)}\big(S_n(R)f_1\oplus S_n(R)f_{m+1}\big)\simeq S_3(R)$.

(2) If $n=2m$, then $1=f_1+\cdots + f_m$. Similarly, we can prove $S_n(R)f_j\simeq S_n(R)f_1$ as left $S_n(R)$-modules for
$2\le j\le m$. Thus $S_n(R)$ is Morita equivalent to $\End_{S_n(R)}\big(S_n(R)f_1\big).$ It is easy to show $\End_{S_n(R)}\big(S_n(R)f_1\big)\simeq f_1S_n(R)f_1\simeq S_2(R)$.

Another proof of (2) can be found in the proof of Theorem \ref{thm2}.

(3) Now, let $n$  be arbitrary. To show that the extension $S_n(R)\subseteq M_n(R)$ is a Frobenius extension, we define
$$\sigma: M_n(R)\lra M_n(R), \quad y\mapsto cyc.$$
Note that $c^2=c$ and $cs=sc$ for all $s\in S_n(R)$. Thus $a +\sigma(a)=a+cac=c(a+cac)c\in S_n(R).$ Visually, if $a=(a_{ij})\in M_n(R)$, then
\[\sigma(a)=\begin{pmatrix}
a_{n n} & a_{n n-1} & \dots & a_{n 2} & a_{n 1} \\
a_{n-1,n} & a_{n-1,n-1} & \dots & a_{n-1,2} & a_{n-1,1} \\
\vdots & \vdots &  & \vdots & \vdots \\
a_{2 n} & a_{2,n-1} & \dots & a_{2 2} & a_{2 1} \\
a_{1 n} & a_{1,n-1} & \dots & a_{1 2} & a_{1 1} \\
\end{pmatrix}_{n\times n}\]

So \[a+\sigma(a)=
\begin{pmatrix}
a_{1 1}+a_{n,n} & a_{1 2}+a_{n,n-1} & \dots & a_{1,n-1}+a_{n 2} & a_{1 n}+a_{n 1} \\
a_{2 1}+a_{n-1,n} & a_{2 2}+a_{n-1,n-1} & \dots & a_{2,n-1}+a_{n-1,2} & a_{2 n}+a_{n-1,1} \\
\vdots & \vdots &  & \vdots & \vdots \\
a_{n-1,1}+a_{2 n} & a_{n-1,2}+a_{2,n-1} & \dots & a_{n-1,n-1}+a_{2 2} & a_{n-1,n}+a_{2 1} \\
a_{n 1}+a_{1 n} & a_{n 2}+a_{1,n-1} & \dots & a_{n,n-1}+a_{1 2} & a_{n,n}+a_{1 1} \\
\end{pmatrix}_{n\times n}\in S_n(R).\]

Now, we define
$$ E: M_n(R)\lra S_n(R),\quad x\mapsto x+\sigma(x)=x + cxc\in S_n(R),$$
$x_i:=e_{i 1}$ and $y_i:=e_{1 i}\in M_n(R)$ for $1\le i\le n$. In the following, we will prove that $(E,x_i,y_i)$ is a Frobenius system, that is, the following (i) and (ii) hold.

(i) $E$ is a homomorphism of $S_n(R)$-$S_n(R)$-bimodules. In fact, we have $E(x+y)=E(x) +E(y).$ Moreover, for $s\in S_n(R)$ and $x\in M_n(R)$, it follows from $sc=cs$ that $E(sx)=sx+csxc = sx+s\cdot cxc =s(x+cxc)=sE(x)$. Similarly, $E(xs)=E(x)s$.

(ii) For any $a\in M_n(R)$, $\sum_i x_iE(y_ia)=a$ and $\sum_i E(ax_i)y_i=a$. Indeed, $ax_i=ae_{i 1}$ is the matrix which has the first column equal to the $i$-th column of $a$ and other columns equal to $0$, and $x_ia=e_{i 1}a$ is the matrix that has the $i$-th row equal to the first row of $a$ and other rows equal to $0$. Also, $y_ia$ is the matrix with the first row equal to the $i$-th row of $a$ and other rows equal to zero, while $ay_i$ is the matrix with the $i$-th column equal to the first column of $a$ and other columns equal to $0$. Thus $x_iE(y_ia)$ is the matrix with the $i$-th row equal to the $i$-th row of $a$ and with other rows equal to $0$. Hence $\sum_i x_iE(y_ia)=a$. Similarly, $\sum_i E(ax_i)y_i=a$.

Thus $(E, x_i,y_i)$ is a Frobenius system and $S_n(R)\subseteq M_n(R)$ is a Frobenius extension. Since $\sum_{i=1}^nx_iy_i=1$, the Frobenius extension is separable by Lemma \ref{lem3}(2)(i). $\square$

\begin{Koro} If $2$ is invertible in $R$, then $S_n(R)\subseteq M_n(R)$ is a split Frobenius extension.

\end{Koro}

{\it Proof.} Since $c^2 =1 \in M_n(R)$, we have $E(\frac{1}{2})= \frac{1}{2} +\frac{1}{2}c^2=1\in S_n(R)$.
Thus the extension is split by Lemma \ref{lem3}(ii).
 $\square$

Next, we show the cellularity of centrosymmetric matrix algebras.

\begin{Theo} \label{thm2} If $R$ is a commutative ring, then

$(1)$ $S_{2m+1}(R)$ is a quasi-hereditary, cellular $R$-algebra with respect to the matrix transpose.

$(2)$ $S_{2m}(R)$ is a cellular $R$-algebra with respect to the matrix transpose.
\end{Theo}

{\it Proof.} Note that $S_n(R)$ is closed under taking transpose and  $f_{i}'=f_{i}\in S_n(R)$ for all $i$.

(1) Let $n=2m+1$ with $m\ge 0$. Then $f_{m+1}S_n(R)f_{m+1}\simeq R$. We consider the ideal $S_n(R)f_{m+1}S_n(R)$ and show that the multiplication map
$$ S_n(R)f_{m+1}\otimes_{R}f_{m+1}S_n(R)\lra S_n(R)f_{m+1}S_n(R)$$
is injective, and therefore an isomorphism. For this purpose, it is sufficient to show that, for any $1\le i,j\le m+1$, the multiplication map
$$ \mu_{ij}: f_iS_n(R)f_{m+1}\otimes_{R}f_{m+1}S_n(R)f_j\lra f_iS_n(R)f_{m+1}S_n(R)f_j$$
is injective. If one of $f_i$ and $f_j$ equals $f_{m+1}$, then $\mu_{ij}$ is obviously injective. So, we assume $1\le i,j\le m$. In this case, $f_iS_n(R)f_{m+1}=Rf_{i, m+1}$, $f_{m+1}S_n(R)f_j=Rf_{m+1, j}$ and
$f_iS_n(R)f_{m+1}S_n(R)f_j=Rf_{i, m+1}f_{m+1, j}=R(f_{ij}+f_{i,n+1-j}).$ Thus the $R$-module $f_iS_n(R)f_{m+1}\otimes_{R}f_{m+1}S_n(R)f_j$ is a free $R$-module of rank $1$. Hence
the multiplication map $\mu_{ij}$ is injective.

Let $J=S_n(R)f_{m+1}S_n(R)$ and $\overline{S_n(R)}:=S_n(R)/J$. Then $J$ has an $R$-basis $\{\alpha_{ij}:=f_{ij}+f_{i,n+1-j}\mid 1\le i,j\le m\}\cup\{\alpha_{m+1,j}:=f_{m+1,j}, \alpha_{j,m+1}:=f_{j,m+1}\mid 1\le j\le m\}$. Let $J'_1=J$ and $J'_2=\{f_{ij}\mid 1\le i,j \le m\}$. Then $S_n(R)$ has an $R$-module decomposition $S_n(R)=J'_1 \oplus J'_2$ such that the matrix transpose preserves $J'_1$ and $J'_2$, respectively. By Lemma \ref{lem2.7}, $J$ is a cell ideal. Now we prove that $\overline{S_n(R)}$ is a cellular algebra with respect to the matrix transpose. It suffices to show $\overline{S_n(R)}\simeq M_m(R)$ as cellular algebras.

Indeed, $\{\bar{f}_{ij}\mid 1\le i,j\le m\}$ is an $R$-basis of $\overline{S_n(R)}$ and $\bar{1}=\bar{f}_1+\cdots +\bar{f}_m$, where $\bar{a}=a+J\in \overline{S_n(R)}$ for $a\in S_n(R)$.
Thus
$$ \overline{S_n(R)}=\begin{pmatrix}
\bar{f}_1\overline{S_n(R)}\bar{f}_1 & \bar{f}_1\overline{S_n(R)}\bar{f}_2 & \cdots & \bar{f}_1\overline{S_n(R)}\bar{f}_m\\
\bar{f}_2\overline{S_n(R)}\bar{f}_1& \bar{f}_2\overline{S_n(R)}\bar{f}_2& \cdots &\bar{f}_2\overline{S_n(R)}\bar{f}_m\\
\vdots & \vdots &\ddots & \vdots\\
\bar{f}_{m}\overline{S_n(R)}\bar{f}_1&\bar{f}_{m}\overline{S_n(R)}\bar{f}_2& \cdots & \bar{f}_m\overline{S_n(R)}\bar{f}_m
\end{pmatrix}.$$
(For basic facts on representing an algebra as a matrix algebra, we refer to, for instance, \cite[Section 1.7]{DK}).
It follows from $f_{i, m+1}f_{m+1, n+1-j}= f_{ij}+ f_{i, n+1-j}$ (see Lemma \ref{lem2}(2)) that $\bar{f}_{ij}=-\bar{f}_{i, n+1-j}$ in $\overline{S_n(R)}.$ Since $f_iS_n(R)f_j=Rf_{ij}+Rf_{i, n+1-j}$, we get $\bar{f}_i\overline{S_n(R)}\bar{f}_j =R\bar{f}_{ij}$ and $\overline{S_n(R)}=\bigoplus_{1\le i,j\le m}R\bar{f}_{ij}.$

If $1\le i,j\le m$, then $\delta_{i,n+1-j}=0$. Otherwise, we would have $2m+2=n+1=i+j\le 2m.$ Thus it follows from Lemma \ref{lem2}(2) that, for $1\le i,j,p,q\le m$,
$$f_{ij}f_{pq}=\delta_{jp}f_{iq}+\delta_{j, n+1-p}f_{i, n+1-q} =\delta_{jp}f_{iq}.$$
Hence $\bar{f}_{ij}\bar{f}_{pq}=\delta_{jp}\bar{f}_{iq}$ in $\overline{S_n(R)}.$  This means $\overline{S_n(R)}\simeq M_m(R)$. Clearly, the transpose on matrices in $\overline{S_n(R)}$ is just the transpose of matrices in $M_m(R)$. Thus this isomorphism is an isomorphism of cellular algebras.

(2) Let $n=2m$ with $m\ge 1$. Then $\delta_{i,n+1-i}=0$ for all $i$, and $\delta_{i,n+1-j}=0$ for all $1\le i, j\le m$. Note that $1=f_1+\cdots + f_m$ with $f_i=e_{ii}+ e_{n+1-i, n+1-i}$. Thus $S_n(R)$ has the matrix representation
$$ S_{n}(R)= \begin{pmatrix}
{f}_1S_n(R)f_1 & f_1S_n(R)f_2 & \cdots & f_1S_n(R)f_m\\
f_2S_n(R)f_1& f_2 S_n(R)f_2 & \cdots &f_2S_n(R)f_m\\
\vdots & \vdots &\ddots & \vdots\\
f_{m}S_n(R)f_1& f_m S_n(R)f_2& \cdots & f_m S_n(R)f_m
\end{pmatrix}.$$
The transpose is then given by $(a_{ij})'=(a'_{ji})$ where $a_{ij}\in f_iS_n(R)f_j$ for all $i, j$.
By Lemma \ref{lem2}(1), we have
$$ f_iS_n(R)f_j=Rf_{ij}\oplus Rf_{i,n+1-j}, \quad f_{i,n+1-i}^2=f_i, \quad 1\le i,j\le m.  $$
Note that the above direct sum of $R$-modules follows from $f_{ij}\ne f_{i, n+1-j}$.
Let $x$ be a generator of the cyclic group $C_2$ of order $2$ and $\varphi$ be the map
$$\varphi: M_m(R[C_2])\lra S_n(R), \quad \big(a_{ij}+b_{ij}x\big)\mapsto \big(a_{ij}f_{ij}+b_{ij}f_{i, n+1-j}\big)$$
where $a_{ij},b_{ij}\in R$ for all $i,j$. Then, by Lemma \ref{lem2}, we can check that $\varphi$ is an algebra homomorphism.  Obviously, $\varphi$ is surjective and injective.
Thus $S_{2m}(R)\simeq M_m(R[C_2])$ as $R$-algebras. Observe that the commutativity of $R$ is not used in the above argument and therefore this isomorphism holds true for any ring $R$.

Note that $\varphi$ commutes with the matrix transpose and that $M_m(R[C_2])$ is a cellular algebra with the transpose as its involution. Thus $S_{2m}(R)$ is a cellular $R$-algebra with respect to the transpose of matrices. $\square$

\medskip
The following corollary describes the algebraic $K$-groups of the centrosymmentric matrix algebras. We refer to \cite{Ros} for basic knowledge on algebraic $K$-theory of rings.

\begin{Koro} $(1)$ Let $K_n(R)$ denote the $n$-th algebraic $K$-group of a ring $R$ in the sense of Quillen. Then  $K_n\big(S_{2m}(R)\big)\simeq K_n(R[C_2])$ and $K_n\big(S_{2m+1}(R)\big)\simeq K_n(R)\oplus K_n(R)$ for $n\ge 0$ and $m\ge 1.$

(2) If $R$ is commutative, then the centre of $S_n(R)$ is $R[c]$ for $n\ge 1$.
\end{Koro}

{\it Proof.} (1) Since Morita equivalences preserve algebraic $K$-groups of rings, $K_n(S_{2m}(R))\simeq K_n(S_2(R))\simeq K_n(R[C_2])$ for all $n\ge 1$ by Theorem \ref{thm1}.

Let $A:=S_{2m+1}(R)$. It follows from the proof of Theorem \ref{thm2} that the ideal $J:=Af_{m+1}A$ is a finitely generated projective $A$-module and $\End_A(_AJ)$ is Morita equivalent to $f_{m+1}Af_{m+1}\simeq R$. Then, by \cite[Corollary 1.5]{CX}, $K_n(A)\simeq K_n(A/J)\oplus K_n(\End_A(_AJ))=K_n(M_m(R))\oplus K_n(R)\simeq K_n(R)\oplus K_n(R)$.

(2) Recall that the center of a ring $S$ is the set $Z(S):=\{x\in S\mid xy=yx, \forall y\in S\}$. Since Morita equivalences preserve the centers of rings, the center of $S_{2m}(R)$ is $R[c]$.
For $S_{2m+1}(R)$, it follows from Hochschild cohomology calculation of quasi-hereditary algebras (see \cite{DX}) that
the centre of $S_{2m+1}(R)$ is a free $R$-module of rank $2$ for $m\ge 1$, thus it is $R[c]$. For $S_1(R)$, we have $c=1$ and $R[c]=R$.
$\square$

Finally, we point out quiver presentations of centrosymmetric matrix algebras over a field.

Let $R$ be a field. The path algebra of a quiver $Q$ over $R$ is denoted by $RQ$. By Theorem \ref{thm1}, up to Morita equivalence, it is sufficient to describe the quivers and relations of $S_i(R)$ for $i=1,2,3.$

Clearly, $S_1(R)=R\;\bullet$. If char$(R)\neq 2$, then $S_2(R)=R(\bullet\quad \bullet)$ and $S_3(R)=M_2(R)\times R\approx R(\bullet\quad \bullet)$.
If char$(R)=2$, then  $$S_2(R)=R\big(\xymatrix{\bullet\ar@(ur,dr)[]^{\alpha}}\big)/(\alpha^2)$$ and
$$S_3(R)=R\big( \xymatrix{
 \bullet\ar@<2.5pt>[r]^{\alpha}& \bullet\ar@<2.5pt>[l]^{\beta}}
\big)/(\alpha\beta).$$

Recall that a finite-dimensional algebra $A$ over a field is representation-finite if there are finitely many non-isomorphic indecomposable left $A$-modules. Thus

\begin{Koro} The centrosymmetric matrix algebra $S_n(R)$ over a field $R$ is representation-finite for all $n\ge 1$.
\end{Koro}

\medskip
{\bf Acknowledgements.} The research work of the authors was partially supported by BNSF(1192004).

\medskip
{\footnotesize

Changchang Xi, School of Mathematical Sciences, Capital Normal University, 100048 Beijing, China

{\tt Email: xicc@cnu.edu.cn}

\medskip
Shujun Yin, School of Mathematical Sciences, Capital Normal University, 100048 Beijing, China

{\tt Email: yinsj9417@163.com}
}
\end{document}